\documentclass[12pt]{article}
\begin{document}
\font\aaa=msbm10 scaled \magstep 0
\begin{center}
{\Large {\bf { The Ricci Flow on Complete Noncompact K\"{a}hler Manifolds}}}
\vskip 0.5cm
 Xi-Ping Zhu
\vskip 0.3cm
Department of Mathematics, Zhongshan University\\
Guangzhou 510275, P. R. China 
\end{center}
\vskip 0.6cm
\begin{abstract}
 In this paper we survey the recent developments of the Ricci flows on
 complete noncompact K\"{a}hler manifolds and their applications in geometry.
\end{abstract}
\vskip 0.6cm
\section{  Introduction}
\vskip 0.3cm The classical uniformization theorem of Riemann
surface states that a simply connected Riemann surface is
biholomorphic to either the Riemann sphere, the complex line or
the open  unit disc. It gives the characterization for the
standard complex structures of one-dimensional K\"{a}hler
manifolds. Unfortunately, a direct analog of this beautiful result
to higher dimension does not exist. For example, it is well-known
that there is a vast variety of biholomorphically distinct complex
structures on $\aaa\mbox{R}^{2n}$ for $n > 1$. This says that the
topological restrictions are no longer distinguishing the standard
complex structures in higher dimensions. Thus, in order to
characterize the standard complex structures for higher
dimensional manifolds, one needs to impose geometric restrictions
on the manifolds. \vskip0.1cm From the point of view of
differential geometry, one consequence of the uniformization
theorem is that a positively curved, compact or noncompact Riemann
surface must be biholomorphic to the Riemann sphere or the complex
line respectively. The higher dimensional version of the
uniformization theorem is to search the similar characterization
for complete K\"{a}hler manifolds with ``positive curvature''.
Such a characterization in the case of compact K\"{a}hler manifold
is the famous Frankel conjecture which says that a compact
K\"{a}hler manifold of positive holomorphic bisectional curvature
is biholomorphic to a complex projective space. This conjecture
was solved by Andreotti-Frankel [11] and Mabuchi [19] in complex
dimensions two and three respectively and the general case was
then solved by Mori [22], and Siu-Yau [33] independently. Thus the
further investigations are naturally led to complete noncompact
K\"{a}hler manifolds of positive holomorphic bisectional
curvature. The following conjecture provides the main impetus.
\vskip 0.3cm \noindent {\bf  Conjecture 1} {\it (Greene-Wu [13],
Yau [37]) A complete noncompact K\"{a}hler manifold of positive
holomorphic bisectional
 curvature is biholomorphic to a complex Euclidean space.}
\vskip 0.3cm Recall from Cheeger-Gromoll-Meyer [3], [14] that a
complete noncompact Riemannian manifold with positive {\it
sectional } curvature is diffeomorphic to a Euclidean space and
recall from Greene-Wu [12] that a complete noncompact K\"{a}hler
manifold with positive {\it sectional } curvature is Stein.
Nonetheless, in the case of positive {\it holomorphic bisectional}
curvature, very little is known about the topology and complex
structure of these manifolds. For example, one does not even know
whether a complete noncompact K\"{a}hler manifold of positive
holomorphic bisectional curvature is simply connected. Moreover it
is also unknown whether a complete noncompact K\"{a}hler  manifold
of positive holomorphic bisectional curvature is Stein, which is a
conjecture of Siu [32]. The Ricci flow introduced by Hamilton [15]
has been found to be an useful tool to understand these manifolds.
\vskip 0.3cm Let $M$ be a complete noncompact K\"{a}hler manifold
with K\"{a}hler metric $g_{\alpha \bar{\beta}}.$ The Ricci flow is
the following evolution equation on the metric
$$
\arraycolsep=1.5pt
\left\{
\begin{array}{rcl}
\displaystyle \frac {\partial g_{_{\alpha \bar{\beta}}}}{\partial t}(x,t) & = &
-R_{_{\alpha \bar{\beta}}}(x,t), \quad x \in M,t>0,
\\[4mm]
\displaystyle g_{_{\alpha \bar{\beta}}}(x,0) & = &g_{_{\alpha \bar{\beta}}}(x), \quad x \in M,
\end{array}
\right.
\eqno (1.1)
$$
where $R_{_{\alpha \bar{\beta}}}(x,t)$ denotes the Ricci curvature
tensor of the metric $g_{_{\alpha \bar{\beta}}}(x,t).$ \vskip 0cm
The first result to the Ricci flow (1.1) on  complete noncompact
K\"{a}hler manifolds is the short time existence obtained by Shi
[28]. The short time existence states that if the curvature of the
initial metric $g_{\alpha \bar{\beta}}$ is bounded, the Ricci flow
(1.1) has a maximal solution $g_{\alpha \bar{\beta}}(\cdot,t)$ on
$[0, t_{\max})$ with $t_{\max} > 0$ and the curvature of
$g_{\alpha \bar{\beta}}(\cdot,t)$ becomes unbounded as $t$ tends
to $t_{\max}$ if $t_{\max} < + \infty.$ It is also known (see [17]
or [31]) that the Ricci flow (1.1) preserves the K\"{a}hlerity of
the metric and that if the initial metric has positive holomorphic
bisectional curvature, the evolving metric still has positive
holomorphic bisectional curvature for all times. Thus to study the
topological and complex structure of a complete noncompact
K\"{a}hler manifold of positive holomorphic bisectional curvature,
we can replace the  K\"{a}hler metric by any one of the evolving
metric of  the Ricci flow (1.1). The most investigations of the
Ricci flow (1.1) on complete noncompact K\"{a}hler manifolds are
mainly devoted to understanding the topology, geometry at infinity
and complex structures of such manifolds with positive holomorphic
bisectional curvature. \vskip 0.1cm In this paper we will describe
the recent developments of the Ricci flow (1.1) on noncompact
K\"{a}hler manifolds and will emphasize on its applications in
geometry. Section 2 discusses the long time existence of the Ricci
flow (1.1). In Section 3 we present the gap phenomena for complete
noncompact K\"{a}hler manifolds of positive holomorphic
bisectional curvature. Further the geometry of these manifolds at
infinity is investigated in Section 4. Finally in Section 5 we
report the recent progresses on the long standing conjecture of
Greene-Wu and Yau. \vskip 0.6cm
\section{  Long Time Existence}
\vskip 0.3cm
Let $(M, g_{\alpha \bar{\beta}})$ be a complete  noncompact K\"{a}hler
manifold with bounded and positive holomorphic bisectional curvature.
Consider the maximal solution $g_{\alpha \bar{\beta}}(\cdot,t)$
$(t \in [0, t_{\max}))$ of the Ricci flow (1.1) with
$g_{\alpha \bar{\beta}}$ as the initial metric. By a direct calculation
from the equation in (1.1), the scalar curvature $R(x,t)$ of the
metric $g_{\alpha \bar{\beta}}(x,t)$ is evolved by
$$\frac {\partial R}{\partial t}(x,t) = \Delta_{g(t)}
R(x,t)+2| R_{\alpha \bar{\beta}}(x,t)|_{g(t)}^2, \quad
\mbox{ on } M \times [0, t_{\max}), \eqno (2.1)$$
where $\Delta_{g(t)}$ is the Laplace-Beltrami operator with respect to
the evolving metric $g_{\alpha \bar{\beta}}(x,t)$. This is a nonlinear
heat equation with superlinear growth. At first sight, one thus believes
that the scalar curvature will generally blow up in finite time.
\vskip 0.1cm
In [31], Shi considered the volume element of the evolving metric
and introduced the following function
$$F(x,t)= \log \frac {\det (g_{\alpha \bar{\beta}}(x,t))}
{\det ( g_{\alpha \bar{\beta}}(x,0))}, \eqno (2.2)$$
on $M \times [0,t_{\max}).$ It can be easily obtained from (1.1) that
$$\frac {\partial F(x,t)}{\partial t}= -R(x,t). \eqno (2.3)$$
Since the holomorphic bisectional curvature of $g_{\alpha \bar{\beta}}
(\cdot,t)$ is positive, it follows that $F(\cdot,t)$ is nonincreasing
in $t$ and $F(\cdot,0)=0.$ And by the equation (1.1) we know that
$$g_{\alpha \bar{\beta}}(\cdot,t) \leq g_{\alpha \bar{\beta}}(\cdot,0),
\quad \mbox{ on } M, \eqno (2.4)$$
for any $t > 0$. Then we have
$$
\arraycolsep=1.5pt
\begin{array}{rl}
e^{F(x,t)} R(x,t)= & \displaystyle g^{\alpha \bar{\beta}}(x,t)
R_{\alpha \bar{\beta}}(x,t) \cdot
\frac {\det (g_{\alpha \bar{\beta}}(x,t))}{\det (g_{\alpha \bar{\beta}}(x,0))}
\\[4mm]
\leq & \displaystyle  g^{\alpha \bar{\beta}}(x,0) R_{\alpha \bar{\beta}}(x,t)
\\[4mm]
= & \displaystyle  - \Delta_{g(0)}F(x,t)+R(x,0).
\end{array}
$$
Combining with (2.3) we obtain
$$e^{F(x,t)} \frac {\partial F(x,t)}{\partial t} \geq  \Delta_{g(0)}
F(x,t)- R(x,0), \eqno (2.5)$$
and
$$\Delta_{g(0)} F(x,t) \leq R(x,0), \eqno (2.6)$$
on $M \times [0, t_{\max}).$
\vskip 0.1cm
In the PDE jargon, if the scalar curvature $R(x,0)$ of the initial metric
satisfies suitable decay conditions, the differential inequalities (2.5)
and (2.6) will give two opposite estimates of $F$ by its average.
Shi [31] observed that the combination of these two opposite estimates give
the following a priori estimate for the function $F$.
\vskip 0.3cm
\noindent {\bf Lemma 2.1 } {\it (Shi [31]) Suppose $(M, g_{\alpha \bar{\beta}})$
is a complete noncompact K\"{a}hler manifold with bounded and
 positive holomorphic bisectional curvature. And suppose there exist positive
 constants $C_1,$ $C_2$ and  $0 < \theta < 2$ such that
$$(i) \quad R(x,0) \leq C_1, x \in M, \hskip 4.5cm$$
$$(ii) \quad \frac 1 {Vol (B_0(x_0,r))} \int_{B_0(x_0,r)}R(x ,0)dx \leq \frac
{C_2}{(1+r)^{\theta}},$$
\vskip 0.05cm
$$\hskip -2cm \mbox{  for all  } x_0 \in M,  0 \leq r < +
\infty,$$
where $B_0 (x_0,r)$ is the geodesic ball of radius $r$ and centered at
$x_0$ with respect to the metric $g_{\alpha \bar{\beta}}(x).$ Then
the function $F(x,t)$ satisfies the estimate
$$F(x,t) \geq -C(t+1)^{\frac {2-\theta}{\theta}}, \quad \mbox{ on }
M \times [0,t_{\max}),$$
where $C$ is a positive constant depending only on
$\theta,$ $C_1$, $C_2$, and the dimension.}
\vskip 0.1cm
The  combination of the above lemma and (2.4) implies
$$g_{\alpha \bar{\beta}}(x,0) \geq g_{\alpha \bar{\beta}}(x,t)
\geq e^{-C(t+1)^{\frac {2-\theta}{\theta}}}g_{\alpha \bar{\beta}}
(x,0), \quad \mbox{ on } M \times [0, t_{\max}). \eqno (2.7)$$
Recall that the Ricci curvature is given by
$$R_{\alpha \bar{\beta}} (x,t) =- \frac {\partial^2}{\partial z_{\alpha}
\partial \bar{z}_{\beta}} \log \det (g_{\alpha \bar{\beta}} (x,t)).$$
Thus the Ricci flow equation (1.1) is the parabolic version of the
complex Monge-Amper$\acute{e}$ equation on the K\"{a}hler manifold. The
inequality (2.7) is corresponding to the second order estimate for the
Monge-Amper$\acute{e}$ equation. It is well-known that the third order
and higher order estimates for the Monge-Amper$\acute{e}$ equation were
developed by Calabi and Yau. Similarly, by adapting the Calabi and
Yau's arguments, Shi proved in [31] that the derivatives and higher order
estimates for $g_{\alpha \bar{\beta}}(x,t)$ are uniformly bounded on
any finite time interval. Thus Shi obtained the following long time existence
result.
\vskip 0.3cm
\noindent {\bf Theorem 2.2} {\it  (Shi [31]) Suppose $(M, g_{\alpha
\bar{\beta}})$ is a complete noncompact K\"{a}hler manifold with
bounded and positive holomorphic bisectional curvature. And suppose there
exist positive constants
$C_2$ and $0 < \theta <2$ such that
$$(ii) \quad \frac 1{Vol(B_0(x_0,r))} \int_{B_0(x_0,r)} R(x,0)dx \leq
\frac {C_2}{(1+r)^{\theta}},$$
$$\hskip -2cm \mbox{ for all } x_0 \in M, 0 \leq r < + \infty.$$
Then the Ricci flow (1.1) has a solution for all $t \in [0, +\infty).$}
\vskip 0.1cm
Later on, Ni-Tam [25] observed that when the scalar curvature of the
initial metric decays faster than linear (i.e., $\theta > 1$ in the
condition (ii)), the above long time existence result can be deduced
rather easily. In this case, one can solve the Poisson equation
$$\Delta_{g(0)} u_0(x) = R(x,0)$$
and then by a Bochner technique of Mok-Siu-Yau [21] the solution
$u(x)$ actually satisfies the Poincar$\acute{e}$-Lelong equation
$$\sqrt{-1} \partial \bar{\partial} u_0(x) = Ric (x,0),$$
where $Ric(x,0)$ is the Ricci form of the initial metric
$g_{\alpha \bar{\beta}}(x).$ Set
$$u(x,t)=u_0 (x)- F(x,t)$$
where $F(x,t)$ is defined in (2.2). Then one can easily check
$$\sqrt{-1} \partial \bar{\partial} u(x,t)= Ric(x,t)$$
and
$$( \frac {\partial}{\partial t}- \Delta_{g(t)})u=0$$
on $M \times [0, t_{\max}).$ Moreover one can get the uniform bound
on the gradient to $u(x,t)$. Inspired by Hamilton [16] (see also Chow
[9]), Ni-Tam [25] looked the quantity
$| \nabla u|_{g(t)}^2+ R(x,t)$ and verified
$$(\Delta_{g(t)}-\frac {\partial}{\partial t})(| \nabla u|_{g(t)}^2+
R(x,t)) \geq 0,$$ which concludes from the maximum principle that
the solution of (1.1) exists for all $t \in [0,+ \infty).$
 \vskip
0.1cm The above long time existence result of Shi suggests us to investigate what kind of decay estimate the curvature of positively curved K\"{a}hler manifolds may possess. The classical Bonnet-Myers theorem says that the Ricci
curvature of a complete noncompact manifold can not be uniformly
bounded from below by a positive constant. In [8] Chen and the
author found that the curvature of a complete noncompact
K\"{a}hler manifold actually satisfies the following decay
estimate.
\vskip 0.2cm \noindent {\bf Theorem 2.3} { \it (Chen-Zhu
[8])  Let $M$ be a complete noncompact K\"{a}hler manifold with
positive holomorphic bisectional curvature. Then for any $x_0 \in
M,$ there exists a positive constant $C$ such that
$$\frac 1{Vol(B(x_0,r))}\int_{B(x_0,r)}R(x) dx \leq \frac {C}{1+r},
\quad \mbox{ for all } 0 \leq r < + \infty,$$ where $R(x)$ is the
scalar curvature of $M$.} \vskip 0.1cm We remark that the constant
$C$ in the above decay estimate depending on the point $x_0 \in
M,$ while the assumption (ii) in Theorem 2.2 needs to be uniform
for all $x_0 \in M$. Nevertheless, these two results motivate us
to make the following conjecture.
 \vskip 0.3cm \noindent {\bf
Conjecture 2}  {\it  Let $(M, g_{\alpha \bar{\beta}})$ be a
complete noncompact K\"{a}hler manifold with positive holomorphic
bisectional curvature. Then the Ricci flow (1.1) with $g_{\alpha
\bar{\beta}}$ as initial metric has a solution for all $t \in [0,
+ \infty).$}
 \vskip 0.3cm In [5], Chen, Tang and the author
obtained the following result to support the above conjecture.
\vskip 0.3cm \noindent {\bf Proposition 2.4 } {\it (Chen-Tang-Zhu
[5]) Let $(M, g_{\alpha \bar{\beta}})$ be a complex
two-dimensional complete noncompact K\"{a}hler manifold with
bounded and positive  holomorphic
 bisectional  curvature. Suppose the volume growth of $M$ is
maximal, i.e.,
$$Vol(B(x_0,r)) \geq c r^4, \quad \mbox{ for all } 0 \leq r < + \infty,$$
for some $x_0 \in M$ and some positive constant $c$. Then the
Ricci flow (1.1) with $g_{\alpha \bar{\beta}}$ as the initial
metric has a solution for all $t \in [0. + \infty).$}
\vskip 0.1cm
The proof of this proposition is an indirect blow-up argument. It
also uses some special features in dimension $2$, such as the
Gauss-Bonnet-Chern formula for the four-dimensional Riemannian
manifolds.
\vskip 0.6cm
\section{ Gap Theorems}
\vskip 0.3cm
From now on we discuss the applications of the Ricci flow (1.1) to the
geometry of complete noncompact K\"{a}hler manifolds. In this section we
are interested in the question that how much the curvature of a complete
noncompact K\"{a}hler manifold of nonnegative holomorphic bisectional
 curvature could have near the infinities. The decay estimate of
Theorem 2.3 tells us that there cannot too much at infinity. At first sight,
there seems to be no restriction in the other direction, no limit on
how less curvature can have. For example, it is easy to construct complete
 K\"{a}hler metrics on $\aaa\mbox{C}$ from real surface of revolution
such that their curvatures are zero outside some compact set,
nonnegative everywhere, and positive somewhere. But it is rather
surprising that the corresponding situation can not occur for
higher dimensions. In [21] Mok, Siu and Yau discovered the
following isometrically embedding theorem. \vskip 0.3cm \noindent
{\bf Theorem 3.1 } {\it (Mok-Siu-Yau [21], Mok [20]) Let $M$ be a
complete noncompact K\"{a}hler manifold of nonnegative holomorphic
bisectional curvature of complex dimension $n \geq 2$. Suppose for
a fixed base point $x_0$,
 $$(i) \qquad Vol(B(x_0,r)) \geq C_1 r^{2n}, \quad 0 \leq r < + \infty,
\hskip2cm $$
 $$(ii) \qquad R(x) \leq \frac {C_2}{1+d(x_0,x)^{2+ \varepsilon}}, \quad x \in M, \hskip 3cm$$
 for some $C_1$, $C_2 > 0$ and for any arbitrarily small positive constant
 $\varepsilon.$ Then $M$ is isometrically biholomorphic to $\aaa\mbox{C}^n$ with
 the standard flat metric. }
 \vskip 0.1cm
 This result shows that there is a gap between the flat metric and
 the other metrics of nonnegative curvature on $\aaa\mbox{C}^n.$
 Their method is to consider the Poincar$\acute{e}$-Lelong
 equation $\sqrt{-1} \partial \bar{\partial} u =Ric.$
 Under the condition (ii) that the curvature has faster than
 quadratic decay, they proved the existence of a bounded solution
 $u$ to the Poincar$\acute{e}$-Lelong equation. By virtue of
 Yau's Liouville theorem on complete manifolds with nonnegative
 Ricci curvature, this bounded plurisubharmonic function $u$ must
 be constant and hence the Ricci curvature must be identically
 zero. This implies that the K\"{a}hler metric is flat because of
 the nonnegativity of the holomorphic bisectional curvature. The
 maximal volume growth assumption (i) in the above Theorem 3.1 is
 made to solve the Poincar$\acute{e}$-Lelong equation. This gap
  result was later generalized in [24] for non-parabolic
  manifolds. In [6], Chen and the author gave a further
  generalization as follows.
  \vskip 0.3cm
  \noindent {\bf Theorem 3.2} {\it  (Chen-Zhu [6]) Suppose $M$ is
  a complete noncompact K\"{a}hler manifold of complex dimension
  $n$ with bounded and nonnegative holomorphic bisectional
  curvature. Suppose there exists a positive function
  $ \varepsilon : \aaa\mbox{R} \to \aaa\mbox{R}$
  with $\lim\limits_{r \to + \infty} \varepsilon(r) =0$, such that
  for any $x_0$,
  $$\frac 1{Vol(B(x_0,r))} \int_{B(x_0,r)} R(x) dx \leq \frac {
  \varepsilon(r)}{r^2}.$$
  Then $M$ is a complete flat K\"{a}hler manifold.}
  \vskip 0.1cm
  The method was to use the Ricci flow. We considered the
  K\"{a}hler metric $g_{\alpha \bar{\beta}}(x)$ in Theorem 3.2 as
  the initial metric and evolved it by the Ricci flow (1.1) to get
  a maximal solution $g_{\alpha \bar{\beta}}(x,t)$ on
  $M \times [0,t_{\max}).$ By using the faster than quadratically
  decay condition in the theorem, we could prove that the solution
  exists for all times and the scalar curvature $R(x,t)$ of the
  evolving metric $g_{\alpha \bar{\beta}}(x,t)$ decays faster
  than linear in time, i.e.,
  $$ \lim\limits_{ t \to + \infty} t R(x,t) =0, \quad
  \mbox{ for } x \in M. \eqno (3.1)$$
  On the other hand, Cao [1] obtained the following Li-Yau  type
  inequality for the scalar curvature of the solution of (1.1)
  with nonnegative bisectional curvature,
  $$ \frac {\partial R}{\partial t}-2 \frac {| \nabla R|^2}{R}+
  \frac {R}{t} \geq 0, \quad \mbox{ on } M \times [0, + \infty),
  \eqno (3.2)$$
  which implies that  the function $tR(x,t)$ is nondecreasing in
  time . Thus the combination of (3.1) and (3.2) concludes that
  $R(x,t) \equiv 0$ on $M \times [0,+\infty).$ Therefore $(M,
g_{\alpha \bar{\beta}})$ must be flat.
\vskip 0.1cm
Quite recently Ni-Tam [26] obtained an improved version of the above
gap theorem which does not require the uniform decay and the boundedness
of the curvature tensor.
\vskip 0.3cm
\noindent {\bf Theorem 3.3 } {\it (Ni-Tam [26])  Let $M$ be a complete
noncompact K\"{a}hler manifold with nonnegative holomorphic bisectional
 curvature. Suppose that $R(x) \leq C(d(x_0,x)^2+1)$ for some
$x_0 \in M$ and $C > 0$, and
$$\int_0^r ( \frac s{Vol(B(x_0,s))} \int_{B(x_0,s)}R(x)dx) ds=o( \log r).$$
 Then $M$ must be flat.}
\vskip 0.6cm
\section{ Volume  Growth and Curvature Decay}
\vskip 0.3cm
We continue to study  the geometric properties of  positively curved
complete noncompact manifolds near the infinities. Let us first consider
the volume growth of the manifolds. When an $m-$dimensional complete
 noncompact Riemannian manifold has nonnegative Ricci curvature, the
classical Bishop volume comparison theorem implies that the volume
growth of geodesic balls is at most as the Euclidean volume growth.
On the other hand, Calabi and Yau [36] showed that the volume
growth of a complete noncompact $m-$dimensional Riemannian manifold with
nonnegative Ricci curvature must be at least of linear, i.e.,
$$Vol (B(x_0,r)) \geq cr, \quad \mbox{ for all } 1 \leq r < + \infty,$$
 where $c$ is some positive constant depending on $x_0 \in M$
and the dimension. In [8], Chen and the author found that the volume
growth estimate for K\"{a}hler manifolds of positive holomorphic
bisectional curvature is at least of half of the real dimension (i.e., the
complex dimension).
\vskip 0.3cm
\noindent {\bf  Theorem 4.1 } {\it (Chen-Zhu [8]) Let $M$ be a complex
$n-$dimensional complete noncompact K\"{a}hler manifold with nonnegative
 holomorphic bisectional curvature. Suppose also its holomorphic
bisectional curvature is positive at least at one point. Then the
volume growth of $M$ satisfies
$$Vol (B(x_0,r)) \geq cr^n, \quad \mbox{ for all } 1 \leq r < +\infty,$$
where $c$ is some positive constant depending on $x_0$ and the
dimension $n$.}
\vskip 0.1cm
Note that Klembeck [18] and Cao [2] presented some complete K\"{a}hler
metrics on $\aaa\mbox{C}^n$ which have positive holomorphic bisectional
curvature everywhere such that the volume of the geodesic ball
$B(O,r)$ centered at the origin $O$ with respect to the K\"{a}hler metric
 grows like $r^n$. Thus the volume growth estimate of Theorem 4.1 is
sharp. We also remark that the assumption that the holomorphic
 bisectional curvature is positive at least at one point is necessary.
Indeed, let $M_1$ be a noncompact convex surface in
$\aaa\mbox{R}^3$ which is asymptotic to a cylinder at infinity.
Clearly $M_1$ is a complete noncompact Riemann surface with
positive curvature and has linear volume growth. And let
$\aaa\mbox{C}P^{n-1}$ be the complex projective space with the
Fubini-Study metric. Then the product $M_1 \times
\aaa\mbox{C}P^{n-1}$ is a complex $n-$dimensional complete
noncompact K\"{a}hler manifold with nonnegative holomorphic
bisectional curvature. But its volume growth is of linear. \vskip
0.1cm By the way, in view of Theorem 4.1, it is naturally raised a
question whether there is a similar volume growth estimate for
Riemannian manifolds with positive sectional curvature. The answer
is negative. More precisely, for each dimension $n \geq 2,$ there
exist $n-$dimensional Riemannian manifolds which have positive
sectional curvature  everywhere but have linear volume growth. In
fact on any bounded, convex and smooth domain $\Omega$ in
$\aaa\mbox{R}^n$, we can choose a strictly convex function $u(x)$
defined over $\Omega$ which tends to $+ \infty$ as $x$ approaches
to the boundary $\partial \Omega$. The graph of the convex
function $u(x)$ is a hypersurface in $\aaa\mbox{R}^{n+1}$, denoted
by $M^n$. Clearly the hypersurface $M^n$ is strictly convex, i.e.,
the second fundamental form is strictly positive definite. It then
follows from the Gauss equation that $M^n$ has strictly
 positive sectional curvature. Since the domain $\Omega$ is bounded,
the volume growth of $M^n$ must be of linear.
 \vskip 0.1cm
 We further consider the  relation between the volume growth
 and the curvature decay. Theorem 3.1 of the previous section due
 to Mok, Siu and Yau  implies that the curvatures of complete
 noncompact K\"{a}hler manifolds with positive holomorphic
 bisectional curvature and maximal volume growth can not decay
 faster than quadratically. On the other hand, Yau predicted in
 [38] that if the volume growth is maximal, then the curvature
 must decay quadratically in certain average sense. By
 combining the decay estimate in Theorem 2.3 and
 Chen-Tang-Zhu [5], we had an affirmative answers to the case
 of  complex two-dimensional K\"{a}hler manifolds.
 Later in [8], Chen and the author gave  a further answer for
 higher dimensions under a more restricted curvature assumption.
 More precisely, we have
 \vskip 0.3cm
 \noindent {\bf Theorem 4.2 } {\it (Chen-Tang-Zhu [5],
 Chen-Zhu [8])  Let $M$ be a complex $n-$dimensional
 complete noncompact K\"{a}hler manifold with
 bounded curvature. Suppose the volume growth of $M$ is maximal,
 i.e.,
 $$ Vol(B(x_0,r)) \geq C_1 r^{2n}, \quad \mbox{ for all } 0
 \leq r < + \infty,$$
 for some $x_0 \in M$ and some positive constant $C_1$.
 Suppose also one of followings holds:
 \vskip 0.2cm
 (i) \quad \hskip 0.2cm if $n=2$, the holomorphic bisectional curvature
 of $M$ is positive;
 \vskip 0.2cm
 (ii) \quad if $n \geq 3$, the curvature operator of $M$ is
 nonnegative.
 \vskip 0.2cm
 \noindent
 Then there exists a constant $C_2 > 0$ such that
 $$\frac 1 {Vol(B(x,r))} \int_{B(x,r)} R(y) dy  \leq C_2
 \frac {\log (2+r)}{r^2}, \quad \mbox{ for all }
 x \in M \mbox{ and } r > 0.$$}
 \vskip 0.2cm
 Our method is not direct working on the K\"{a}hler metric.
 Instead we use the K\"{a}hler metric as initial data and
 evolve it by Ricci  flow (1.1). By Theorem 2.3 we know that
 the curvature decays to zero at infinity in average sense. From the
 evolution equation of (1.1) we see that the deformation of the
 metric at infinity is very small. In particular the maximal
 volume growth condition of the initial metric is preserved under
 the evolution of the Ricci flow (1.1), i.e.,
 $$Vol_t(B_t(x,r)) \geq C_1 r^{2n}, \quad \mbox{ for all }
 r > 0, x \in M \mbox{ and } t \in [0, t_{\max}),  \eqno (4.1)$$
 with the same constant $C_1$ as in the initial metric, where
 $B_t(x,r)$ is the geodesic ball of radius $r$ centered at $x$
 with respect to the evolving metric $g_{\alpha \bar{\beta}}(
 \cdot,t)$, and the volume $Vol_t$ is also taken with respect
 to the metric $g_{\alpha \bar{\beta}}(\cdot,t).$ By combining
 with the local injectivity radius estimate of
 Cheeger-Gromov-Taylor [4] we deduce
 $$ inj (M,g_{\alpha \bar{\beta}}(\cdot,t)) \geq \frac {\beta}
 {\sqrt{R_{\max}(t)}}, \quad \mbox{ for } t \in [0, t_{\max}),
 \eqno (4.2)$$
 for some positive constant $\beta > 0$, where $R_{\max} (t)
 =\sup \{ R(x,t) | x \in M \},$
 $ t \in [0, t_{\max}).$
 \vskip 0.1cm
 After obtaining the above injectivity estimate, we can use
 rescaling arguments to analysis the asymptotic behavior of the
 solution $g_{\alpha \bar{\beta}}(\cdot,t)$ near the
 maximal time $t_{\max}$. According to Hamilton [17], the maximal
 solution $g_{\alpha \bar{\beta}}(\cdot,t),$
 $t \in [0,t_{\max}),$ of (1.1) is of either one of the following
 types:
 $$
\arraycolsep=1.5pt
\begin{array}{rl}
\mbox{Type I:} & \quad t_{\max} < + \infty \mbox{ and }
\sup(t_{\max}-t) R_{\max}(t) < +\infty;
\\[4mm]
\mbox{Type II:} & \quad \mbox{either } t_{\max} < +\infty \mbox{
and } \sup (t_{\max}-t) R_{\max}(t)= +\infty,
\\[4mm]
& \quad \mbox{or } t_{\max}= +\infty \mbox{ and } \sup t R_{\max}
(t) = + \infty;
\\[4mm]
\mbox{Type III:} & \quad t_{\max}= + \infty \mbox{ and } \sup t
R_{\max} (t) < + \infty.
\end{array}
$$
Based on the injectivity radius estimate (4.2) we first blow up
the maximal solution into a limiting model and we then blow down
the limiting model at  the infinity to get a rescaling limit which
is splitted as the product of a nonflat solution of the Ricci flow
with some Euclidean space. Since the maximal volume growth
condition (4.1) is invariant under the rescaling, the nonflat
factor must also have maximal volume growth. Thus we can repeat
these rescaling arguments until deriving a complex one-dimensional
nonflat  solution of the Ricci flow with maximal volume growth. If
the original maximal solution were of Type I or II, the rescaling
limit should exist for all $ t \in ( - \infty, 0].$ But only
complex one-dimensional nonflat solution of the Ricci flow of
nonnegative curvature on $(- \infty,0]$ is either the standard
round sphere or the cigar soliton. Since both the round sphere and
the cigar soliton are not of maximal volume growth, we conclude
that the maximal solution must be of type III, i.e., the solution
exists for all $t \in [0,+\infty)$ and satisfies
$$ 0 \leq R(x,t) \leq \frac {C}{1+t}, \quad \mbox{ on }
M \times [0,+ \infty), \eqno (4.3)$$ for some positive constant
$C$. \vskip 0.1cm The scalar  curvature $R(x,t)$ of the evolving
metric $g_{\alpha \bar{\beta}}(x,t)$ satisfies the nonlinear heat
equation (2.1). Intuitively, the Harnack inequality of heat
equation bridges the time decay with the space decay. By using the
time decay estimate (4.3) on the evolving curvature we are indeed
able to prove that the curvature of the initial metric satisfies
the quadratically decay estimate of Theorem 4.2.
 \vskip 0.6cm
\section{ Uniformization  Theorems}
\vskip 0.3cm
 The central question in the study of positively curved complete
 noncompact K\"{a}hler manifolds is the conjecture of Greene-Wu
 and Yau stated in Section 1. A weaker version of this conjecture
 is the following form.
 \vskip 0.3cm
 \noindent {\bf Conjecture 1'} {\it  (Greene-Wu [13], Yau [37])
  A complete noncompact K\"{a}hler manifold of positive sectional
  curvature is biholomorphic to a complex Euclidean space.}
  \vskip 0.1cm
  Since the topology and differential structure of noncompact
  manifolds with positive {\it sectional} curvature is well understood
  from Cheeger-Gromoll-Meyer [3], [14], the above weaker version
  is more concentrative on the complex structure.
  \vskip 0cm
  The first partial affirmative answer to Conjecture 1' was given
  by Mok [20]  for complex two-dimensional manifolds with maximal
  volume growth in the following theorem.
  \vskip 0.3cm
  \noindent {\bf Theorem 5.1 } {\it (Mok [20])  Let $M$ be a
  complete noncompact K\"{a}hler manifold of complex  dimension
  $n$ with positive holomorphic bisectional curvature. Suppose
  there exist positive  constants $C_1$ and $C_2$ such that
  for a fixed base point $x_0$,
  $$(i) \quad Vol(B(x_0,r)) \geq C_1 r^{2n}, \quad 0 \leq r < +
  \infty,$$
  $$ (ii) \quad R(x) \leq \frac {C_2}{1+ d(x_0,x)^2}, \quad \mbox{
  on } M. \hskip 1.2cm$$
  Then $M$ is biholomorphic to an affine algebraic variety.
  Moreover, if in addition the complex dimension $n=2$ and
  $$(iii) \quad \mbox{the sectional curvature of}\quad M\quad \mbox{is positive},$$
  then  $M$ is holomorphic to $\aaa\mbox{C}^2$.}
  \vskip 0.3cm
  Note that  the conditions (i) and (ii) imply the integral
  $\int_M {Ric}^n,$ the analytic Chern number $c_1(M)^n$, is
  finite. The following result of To [35] gave a generalization of
  the above Mok's result to nonmaximal volume growth manifolds.
  It was proved in [35] that if $M$ is a complete noncompact
  K\"{a}hler manifold of positive holomorphic bisectional
  curvature and suppose for some base point $x_0 \in M$ that
  there exist positive $C_1$, $C_2$ and $p$ such that
  $$(i)' \quad \int_{B(x_0,r)} \frac 1{(1+d(x_0,x))^{np}} dx
  \leq C_1 \log (r+2), \quad 0 \leq r < + \infty,$$
  $$(ii)' \quad R(x)  \leq \frac {C_2}{1+d(x,x_0)^p}, \quad
  \mbox{on } M \hskip 5cm $$
  $$(iv) \quad c_1(M)^n = \int_M {Ric}^n < + \infty, \hskip 5.6cm$$
  then $M$ is quasi-projective. Moreover, if in addition the
  complex  dimension $n=2$ and
  $$(iii)\quad the \quad sectional \quad curvature \quad of \quad M \quad is \quad positive,$$
  then $M$ is biholomorphic to
  $\aaa\mbox{C}^2.$
  \vskip 0.1cm
  It is likely that the assumption (iv) is automatically satisfied
  for complete K\"{a}hler manifolds with positive sectional
  curvature. At least in the complex two-dimensional case, there
  holds the generalized Cohn-Vossen  inequality
  $$c_2(M) = \int_M {\Theta} \leq {\bf \cal{\chi}}( \aaa\mbox{R}^4)
   <+\infty$$
  where $\Theta$ is the Gauss-Bonnet-Chern integrand. In view of
  Miyaoka-Yau inequality on the Chern numbers, it is reasonable to
  believe that one can get the finiteness of $c_1(M)^2$ from that
  of $c_2(M).$ Meanwhile in views of Demailly's holomorphic Morse inequality
  [10] and the $L^2$-Riemann-Roch inequality of Nadel-Tsuji [23]
  (see also Tian [34]), the assumption (iv) is a natural condition
  for a complete K\"{a}hler  manifold to be a quasi-projective
  manifold. However the assumptions on the curvature decay and the
  volume growth are more problematic since they demand the
  geometry of the K\"{a}hler manifold at infinity to be somewhat
  uniform. In the recent work [7], Chen and the author showed that
  the assumption (iv) alone is sufficient to give an affirmative
  answer to Conjecture 1' for complex two-dimensional case.
  \vskip 0.3cm
  \noindent {\bf Theorem 5.2} {\it  (Chen-Zhu [7]) Let $M$ be a
  complex $n-$dimensional complete noncompact K\"{a}hler manifold
  with bounded and positive sectional curvature. Suppose
  $$ c_1(M)^n = \int_M {Ric}^n < + \infty.$$
  Then $M$ is biholomorphic to a quasi-projective variety, and in
  case of complex dimension $n=2$, $M$ is biholomorphic to
  $\aaa\mbox{C}^2$.}
  \vskip 0.2cm
  We now come to Conjecture 1. Let us recall a
  remark in [6]. In [29] and [30], Shi used the Ricci flow to
  construct a flat K\"{a}hler metric on $M$ and then claimed
  that $M$ is biholomorphic to $\aaa\mbox{C}^n$. It was point out
  in [6] that from [29] it is unclear why the property of
  completeness is true. Thus one can only get biholomorphic
  embedding of $M$ into a domain of $\aaa\mbox{C}^n$.
  \vskip 0.1cm
  To the knowledge of the author, the only affirmative answer to
  Conjecture 1 on the level of holomorphic bisectional curvature
  is the following result.
  \vskip 0.3cm
  \noindent {\bf Theorem 5.3} { \it (Chen-Tang-Zhu [5], Chen-Zhu [8])
Let $M$ be a complete noncompact  complex two-dimensional K\"{a}hler
manifold of bounded and positive holomorphic bisectional curvature.
Suppose the volume growth of $M$ is maximal, i.e.,
$$Vol(B(x_0,r)) \geq c r^4, \quad \mbox{ for all } 0 \leq r < +\infty,$$
for some point $x_0 \in M$ and some positive constant $c$. Then $M$ is biholomorphic to the complex Euclidean space $\aaa\mbox{C}^2$.}
\vskip 0.2cm
We discribe the proof  of this theorem as follows. The argument there was
 divided into three parts. In the first part, we showed that $M$ is a Stein
manifold homeomorphic to $\aaa\mbox{R}^4$. For this, we evolved the
K\"{a}hler metric  on $M$ by the Ricci flow (1.1). By noting that the
underlying complex structure, K\"{a}hlerness and positivity of holomorphic
bisectional curvature are  preserved under the Ricci flow, the K\"{a}hler
metric in the theorem can be replaced by any one of the evolving
metric. The advantage is that the properties of the evolving metric are
improving during the flow. As seen in the proof of Theorem 4.2, the
maximal volume growth condition is uncharged and the curvature of
evolving metric  decays linearly in time under the Ricci flow (1.1).
These imply that the injectivity radius of the evolving metric is getting
bigger and bigger and any geodesic ball with radius less than half of the
injectivity radius is almost pseudoconvex. By a perturbation argument
 we were then able to modify these geodesic balls to a sequence of exhausting pseudoconvex domains of $M$ such that any two of them form a
Runge pair. From this, it follows readily that $M$ is a Stein manifold
homeomorphic to $\aaa\mbox{R}^4$.
\vskip 0.2cm
In the second part of the proof, we considered the algebra $P(M)$ of
holomorphic functions of  polynomial growth on $M$ and we would prove that
 its quotient field has transcendental degree two over $\aaa\mbox{C}$.
For this, we first needed to construct two algebraically independent holomorphic  functions in the algebra $P(M)$. Using the $L^2$ estimate
of $\bar{\partial}$ operator, it suffices to construct a strictly plurisubharmonic function of logarithmic growth on $M$. As already known
in Mok-Siu-Yau [21], such a strictly plurisubharmonic function of logarithmic growth can be obtain by  solving the Poincar$\acute{e}$-Lelong
equation when the curvature decays in space at least  quadratically.
Fortunately in Theorem 4.2 we had deduced that the curvature of the
initial metric decays quadratically in space in the average sense and
this turns out to be  enough to insure the existence of a strictly
plurisubharmonic function of logarithmic growth. Next , by using
the time decay estimate and the injectivity radius estimate of the
evolving metric, we proved that the dimension of the space of holomorphic
functions in $P(M)$ of degree at most $p$ is bounded by a  constant times
$p^2$. Combining this  with the existence of two algebraically independent
 holomorphic functions in $P(M)$ as above, we could prove that the
quotient field $R(M)$ of $P(M)$ has transcendental degree over
$\aaa\mbox{C}$ by a classical argument of Poincar$\acute{e}$-Siegel.
In other words, $R(M)$ is a finite extension field of some
$\aaa\mbox{C}(f_1,f_2),$ where $f_1,$ $f_2 \in P(M)$ are algebraically
independent over $\aaa\mbox{C}$. Then from the primitive element theorem,
we deduced $R(M) = \aaa\mbox{C} (f_1, f_2, g/h)$ for some $g, $ $h \in
P(M)$. Hence the mapping $F: M \to \aaa\mbox{C}^4$ given by
$F=(f_1, f_2, g,h)$ defines a birational map from $M$ into some irreducible
affine algebraic subvariety $Z$ of $\aaa\mbox{C}^4.$
\vskip 0.2cm
In the last part of proof, we basically followed the approach of  Mok in [20] to establish a holomorphic map from $M$  onto a
quasi-projective variety by desingularizing the map $F$. Our essential
contribution in this part was to establish uniform estimates on the
multiplicity and the number of irreducible components of the zero divisor
of a holomorphic function in $P(M)$. Again the time decay estimate of the
Ricci flow played a crucial role in the arguments. Based on these estimates,
we could show that the mapping $F: M \to Z$ is almost surjective in the
sense that it can miss only a finite number of subvarieties in $Z$, and
can be desingularized by adjoining a finite number of holomorphic
functions of polynomial growth. This shows that $M$ is a quasi-projective
variety. Finally, by combining with the fact that $M$ is homeomorphic to
$\aaa\mbox{R}^4$, we conclude that $M$ is indeed biholomorphic to
$\aaa\mbox{C}^2$ by a theorem of Ramanujam [27] on algebraic surfaces.

\end{document}